\documentclass{amsart}
\usepackage{amsfonts,amssymb,graphicx}
\newtheorem{theorem}{Theorem}
\newtheorem{lemma}[theorem]{Lemma}
\newtheorem{corollary}[theorem]{Corollary}

\newcommand{\RR}{{\mathbb{R}}}
\begin{document}
\title{The size of spanning disks for polygonal curves}
\author{Joel Hass}
\address{Department of Mathematics,
University of California, Davis California 95616}
\email{hass@math.ucdavis.edu}
\thanks{Partially supported by NSF grant DMS-0072348 and
the Institute for Advanced Study.}
\author{Jack Snoeyink}
\thanks{Partially supported by grants from NSERC}
\address{Department of Computer Science, University of North Carolina, Chapel Hill, NC 27514}
\email{snoeyink@cs.unc.edu}
\author{William P. Thurston}
\address{Department of Mathematics, University of California, Davis
California 95616}
\email{wpt@math.ucdavis.edu}
\thanks{Partially supported by NSF grant DMS-0072540.}
\date{\today}

\begin{abstract}
For each integer $n \ge 0$, there is a closed, unknotted, polygonal
curve $K_n$ in $\RR^3$ having less than
$10n+9$ edges, with the property that any Piecewise-Linear triangulated
disk spanning the curve contains at least $2^{n-1}$ triangles.
\end{abstract}

\maketitle

\section{Introduction.}
Let $K$ be a closed polygonal curve in $\RR^3$ consisting of $n$ line segments.
Assume that $K$ is unknotted, so that it is the boundary of an embedded
disk in $\RR^3$. This paper considers the question: How many triangles
are needed to triangulate a Piecewise-Linear (PL) spanning
disk of $K$?  
The main result, Theorem~\ref{exp.disks} below,
exhibits a family of unknotted polygons with $n$ edges, $n \to \infty$,
such that the minimal
number of triangles needed in any triangulated spanning disk
grows exponentially with $n$.  
More specifically, we
construct a sequence of unknotted simple closed curves $K_n$
in $\RR^3$ having the following properties for each $n \ge 0$.
\begin{itemize}
\item The curve $K_n$ is an unknotted polygon with at most $10n+9$ edges. 

\item Any PL embedding of a triangulated disk into $\RR^3$
with boundary $K_n$ contains at least $2^{n-1}$ triangular faces.
\end{itemize}

\begin{figure}[hbtp]
\centering
\includegraphics[width=.4\textwidth]{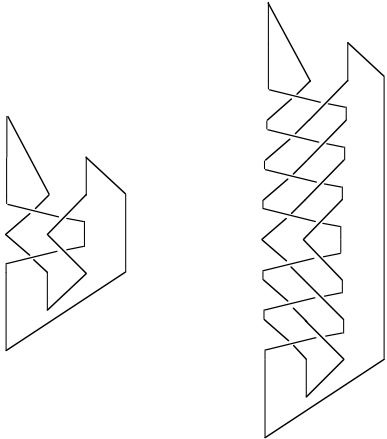}
\caption{The curves $K_1$ and $K_3$.}
\label{Kn}
\end{figure}

The polygons  $K_1$ and $K_3$ are pictured in Figure~\ref{Kn}.

The existence of these curves
is related to the complexity of certain topological algorithms.
Algorithms to test knot triviality by a
search for embedded PL spanning disks
are searching for disks that can be
exponentially more complicated than their
boundary curves. Algorithms of this type include
those described in
\cite{Hak61},\cite{BH97},\cite{HL},\cite{HLP2},\cite{Galatolo}.
Some approaches to problems in computational group theory, such
as the word problem,
are also based on a search for a spanning disk, and may face
similar difficulties.

The lower bound given in our examples
can be compared with the following upper
bound: the results of \cite{HLP2} and \cite{HLT}
show that every unknotted polygon
with at most $n$ edges in $\RR^3$ bounds a PL embedded
triangulated disk which has at most
$C^{n^2}$ triangles, where $C > 1$ is a constant
independent of $n$. The exponent $n^2$ comes from the
requirement that the polygon be embedded in the 1-skeleton of a
triangulated $3$-manifold. A triangulation that contains
an $n$-edge polygon
and using $O(n^2)$ tetrahedra
always exists,
and this bound cannot always be improved, see Avis and El Gindy \cite{AE:87}.
On the other hand, \cite{AE:87} also shows that
a set of points in general position, i.e. one for which
no four points lie on a plane and no three on a line, can be
triangulated using $O(n)$ simplices.
It seems plausible  that one could obtain an improved
upper bound of $C^{n}$ triangles for a PL spanning disk
of a polygon whose vertices are in general position.

A result similar to the one proved in this paper was announced in
\cite{Sno90}, but the geometric analysis suggested there
seems difficult to establish rigorously.
We consider here a different set of polygonal
curves $K_n$ than those used in \cite{Sno90},
and to establish their properties use topological arguments based on
ideas from the classification of diffeomorphisms of surfaces
\cite{Thurston} and from Morse theory~\cite{Milnor}.

Although the main result concerns PL curves and surfaces,
in some parts of the proof it becomes convenient to
work with smooth surfaces and smooth mappings.
This allows use of basic results from smooth Morse theory.
The arguments could be carried out
entirely in the PL context, at the expense of using less well known
versions of Morse theory.
Passing between the PL and smooth settings
is achieved by approximating
PL maps by smooth maps.

\section{Construction of $K_n$} 

\begin{figure}[hbtp]
\centering
\includegraphics[width=.6\textwidth]{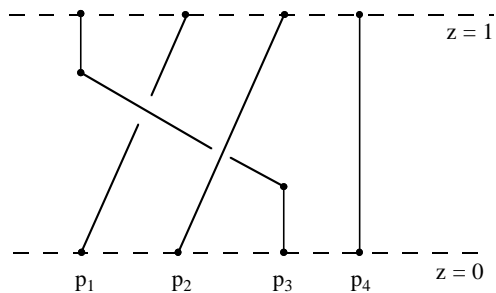}
\caption{The braid $\alpha$.
The four arcs trace out the motion of $p_1, p_2, p_3, p_4$
under a continuous family of plane homeomorphisms
beginning with the identity at $z=1$ and ending with $\varphi$ at $z=0$.
All the arcs lie in the $xz$-plane except the one running from $p_1$ to
$p_2$, which has two vertices off this plane.}
\label{alpha}
\end{figure}

\begin{figure}[hbtp]
\centering
\includegraphics[width=.6\textwidth]{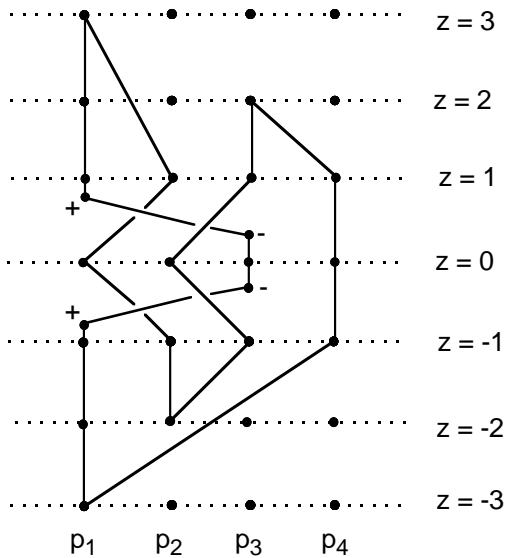}
\caption{A closer look at the unknot $K_1$. All segments lie in the $xz$-plane except
three segments from $p_1$ to $p_3$ between $z=1$ and $z=0$ and
three segments from  $p_3$ to $p_1$ between $z=0$ and $z=-1$.
These have $y$ coordinates above and below the $xz$-plane,
as indicated by the $+$ and $-$ symbols.}
\label{K1}
\end{figure}

We now describe how to construct the unknotted curves $K_n$.
The curve $K_0$
is contained in the $xz$-plane, see Figure~\ref{K1}.
The construction of $K_n$  begins with the PL 4-braid
$\alpha$ depicted in Figure~\ref{alpha}, where
$\alpha = \sigma_1 \sigma_2 ^{-1}$ in terms of the standard generators 
$ \sigma_1,  \sigma_2,  \sigma_3$ of the braid group on four strands, see \cite{Birman}.
This braid consists of four arcs running between the planes
$\{z=1\}$ and $\{z=0\}$, along each of which $z$ is monotonically
decreasing.
The planes $\{z=0\}, \{z=1\}$ each intersect $\alpha$
at four points.
We arrange these points along the $x$-axis and label them by
$p_1=(-2,0),\ p_2=(-1,0),\ p_3=(1,0),\ p_4=(2,0) $.
In this labeling we only consider the $xy$-coordinates.

\begin{figure}[hbtp]
\centering
\includegraphics[width=.8\textwidth]{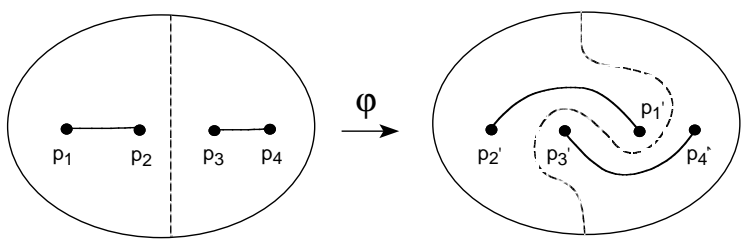}
\caption{The diffeomorphism $\varphi$ associated to $\alpha$.
The diffeomorphism
is the identity outside a disk containing $p_1,p_2,p_3,p_4$.
The locations of the marked points remains fixed, $p_2' = p_1, p_3' = p_2$.}
\label{phi}
\end{figure}

A diffeomorphism $\varphi$ of the 4-punctured plane
$\RR^2 \backslash \{p_1,p_2,p_3,p_4\}$
is associated to the braid $\alpha$.
This diffeomorphism is induced by taking the punctured plane
at level $z=1$ and sliding it down the braid to level $z=0$.
Its action on the plane in indicated in Figure~\ref{phi}.
The action is the identity outside a disk of radius three around the origin.

The curve $K_n$ is formed from an iterated braid
$$
\beta_n = \alpha^{-n} \circ \alpha^{n},
$$
running between the planes $\{z=n\}, \{z=-n\}$.
Between each pair of planes $\{z=k\}$
and $\{z=k+1\}$, $K_n$ consists of a single copy of
$\alpha$ for $0 \le k \le n-1$ and a single copy of
$\alpha^{-1}$ for $-n \le k \le -1$.
In the braid group, $\beta_n$ is equivalent to the trivial
4-braid, which consists of four parallel
vertical segments.
The construction of $K_n$ is completed by
appropriately connecting together the four strands at
the upper and lower ends to form a closed curve,
as shown in Figure~\ref{K1}.
Above the plane  $\{z=n\}$ we
add a pair of line segments from $p_1$ at height
$z=n+2$ to each of $p_1$ and $p_2$ at height
$z=n$, and from  $p_3$ at height
$z=n+1$ to each of $p_3$ and $p_4$ at height 
$z=n$.
Similarly, below the plane
$\{z= -n\}$ we add a pair of line segments from
$p_2$ at height
$z=-n-1$ to each of $p_2$ and $p_3$ at height 
$z=-n$, and from  $p_2$ at height
$z=-n-2$ to each of $p_1$ and $p_4$ at height $z=-n$.
Because the braids $\alpha^n$ and $\alpha^{-n}$ cancel in the
braid group, it is clear that $K_n$ is unknotted.

Our main result is the following:

\begin{theorem}
\label{exp.disks}
For each $n \ge 0$,
\begin{enumerate}
\item $K_n$ is unknotted.
\item $K_n$ contains at most $10n+9$ edges.
\item Any piecewise-smooth embedded disk spanning $K_n$ intersects
the $y$-axis in at least $2^{n-1}$ points.
\item Any embedded PL triangulated disk $D_n$ bounded by $K_n$
contains at least $2^{n-1}$ triangles.
\end{enumerate}
\end{theorem}

The condition (3) that the disk intersects a
line many times implies condition (4),
that it contains many triangles, since each triangle
can intersect a line transversely at most once.

\begin{figure}[hbtp]
\centering
\includegraphics[width=.6\textwidth]{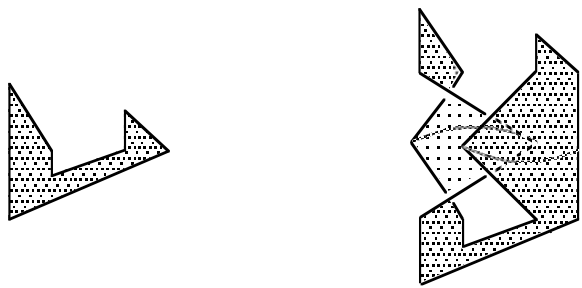}
\caption{The standard disks $F_0$ and $F_1$.}
\label{standard}
\end{figure}

We prove Theorem~\ref{exp.disks} in \S\ref{std}. We
first construct a {\em standard spanning disk} for
$K_n$, which we call $F_n$.
Figure~\ref{standard} shows $F_0$ and $F_1$.
To understand the behavior of $F_n$ and other disks spanning $K_n$,
we prove some facts about diffeomorphisms and ``train tracks''.
These are applied in \S\ref{invt} to count the
intersections of the $y$-axis and $F_n$.
In \S\ref{comb} we use Morse Theory to show that any other
spanning disk is at least as
complicated, along the $y$-axis, as the standard disk.

\section{Construction of a standard spanning disk} 
\label{std}

\begin{figure}[hbtp]
\centering
\includegraphics[width=.8\textwidth]{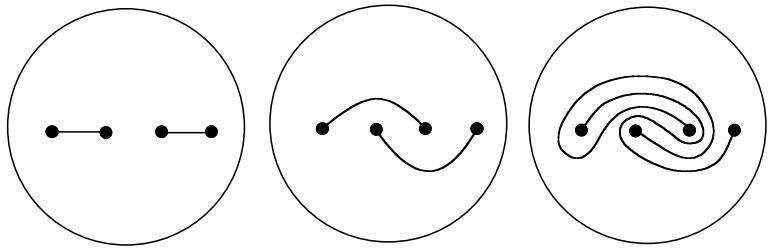}
\caption{Arcs of intersection of a standard disk
with the planes $\{z=n\}, \{z=n-1\}$ and $\{z=n-2\}$.}
\label{arcs}
\end{figure}

In this section we describe how to
construct for each $K_n$ a particular smooth spanning disk $F_n$.
This {\em standard disk} intersects each plane $\{ z=c\}, \ -n-1 < c < n+1$
in two arcs, which are embedded and disjoint.
At $z= \pm n$ the arcs lie along the $x$-axis, joining
$p_1, p_2$ and $p_3,p_4$ respectively.
For $n \ge 2$ these arcs are shown in Figure~\ref{arcs} at 
heights $z=n, z=n-1$ and $z=n-2$.
In  Figure~\ref{arcs} the four arcs appear in the three pictures
in order (1,2,3,4), (2,3,1,4) and (3,1,2,4), read left to right.

Above the plane $\{z=n\}$ the standard disk 
consists of two triangles
in the $xz$-plane, one with a base along the segment
from $p_1$ to $p_2$ and one with a base along the segment
from $p_3$ to $p_4$.
Below $z=-n$ it is bounded by a six-sided polygon
in the $xz$-plane meeting
$\{z=n\}$ along two segments, one running from 
$p_1$ to $p_2$ and one from $p_3$ to $p_4$.
Between $\{z=n\}$ and $\{z=-n\}$ the standard disk
twists so that its boundary follows $K_n$,
as made precise below.

The arcs in the first disk of Figure~\ref{arcs}
are taken by $\varphi$ to the arcs in the second disk in the Figure, and those
in turn are taken by $\varphi$ to the arcs in the rightmost disk. 
The arcs of the braid indicate the motion of the disk
in the process of sliding from $z=1$ to $z=0$.
A composition of a counterclockwise half-twist
interchanging the first two punctures, followed by a
clockwise half-twist interchanging the second two punctures gives $\varphi$. 
Corresponding to the braid $\alpha^n$ is the iterate $\varphi^n$ of $\varphi$.

\begin{figure}[hbtp]
\centering
\includegraphics[width=.1\textwidth]{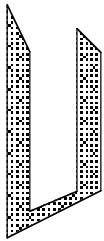}
\caption{A planar disk $L_n$, lying between
$\{ z=n+2 \}$  and $\{ z = -n-2 \}$, is twisted to form $F_n$.}
\label{Ln}
\end{figure}

We now give a precise description of the
construction of $F_n$, based on $\varphi$.
Begin with a planar polygonal curve bounding
a disk $L_n$ in the $xz$-plane, formed as follows:
Take vertical segments from 
$(-2, 0, -n)$ to $(-2, 0, n)$,
$(-1, 0, -n)$ to $(-1, 0, n)$,
$(1, 0, -n)$ to $(1, 0, n)$ and
$(2, 0, -n)$ to $(2, 0, n)$.
At the top, add a line segment from $(-2, 0, n+2)$
to each of $(-2, 0, n)$, $(-1, 0, n)$ and 
from $(1, 0, n+1)$ to each of $(1, 0, n)$, $(2, 0, n)$.
At the lower end. add a line segment connecting
$(-2, 0, -n-2)$ to each of $(-2, 0, -n)$,
$(2, 0, -n)$ and $(-1, 0, -n-1)$ to each of $(-1, 0, -n)$,
$(1, 0, -n)$, as shown in Figure~\ref{Ln}.

The standard disk $F_n$ is the image of this planar disk
$L_n$ under a diffeomorphism $J_n: \RR^3 \to \RR^3$, that
preserves the $z$-coordinates of points, and
carries the boundary of $L_n$ to $K_n$.
The diffeomorphism $\varphi$ 
is isotopic to the identity map on the plane.
So there is a continuous
family of diffeomorphisms of the plane $\varphi_t, \ 0 \le t \le 1$, with 
$\varphi_1 =$ identity and $\varphi_0 = \varphi$. 
Define a diffeomorphism $G:\RR^2 \times [0,1] \to \RR^2 \times [0,1]$
by $G(x,y,t) = (\varphi_t(x,y), t)$.  Then
$G$ carries the vertical line segments in $\RR^2 \times [0,1]$ with
$xy$-coordinates $(-2,0), (-1,0), (1,0),(2,0)$
to the braid $\alpha$.
Also define a diffeomorphism
$H:\RR^2 \times [0,1] \to \RR^2 \times [0,1]$,
by $H(x,y,t) = (\varphi(x,y), t)$.  
This extends $\varphi$ to $\RR^2 \times [0,1]$
as a product.
The diffeomorphism $J_n$ is defined to be the identity for
$z \ge n$ and $z \le -n$. For $0 \le k \le t \le k+1 \le n$,
$J_n(x,y,t) = G \circ H^{n-k-1}(x,y,t-k) + (0,0,k)$,
and for $-n \le k-1 \le t \le k \le 0$,
$J_n(x,y,t) = G \circ H^{n+k-1}(x,y,k-t) + (0,0,2t-k)$.  

A Morse function $f$ on a smooth closed manifold has a finite number of critical
points $\{c_i]$ with distinct values under $f$. 
A Morse function on a manifold with boundary is a Morse function 
when restricted to both the
boundary and the interior of the manifold.  Critical points are either
interior critical points or boundary critical points.
A Morse function on a disk has at least two critical points, one maximum and one minimum,
and if there are exactly two critical points then both must occur on the boundary,
since there is a maximum and minimum value for the restriction to the boundary.
The construction of $F_n$ gives $z$ as a Morse function on $F_n$
that has four critical points.
Two are maxima, at $z=n+2$ and $z=n+1$,
one is a minimum, at $z=-n-2$, and one is a saddle point,
at $z=-n-1$. All four critical points lie on the boundary of $F_n$.

\section{An invariant train track for $\varphi$} 
\label{invt}

\begin{figure}[hbtp]
\centering
\includegraphics[width=.8\textwidth]{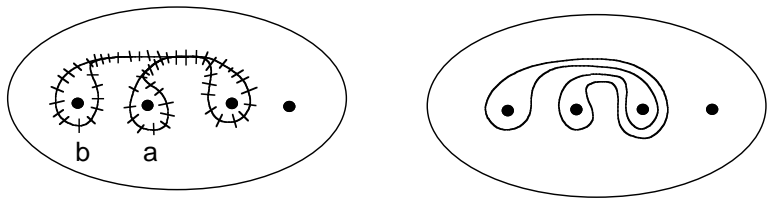}
\caption{A train track $T$, showing its fibered neighborhood, and a
curve that it carries with weights $a=1$ and $b=1$. Other weights
are determined by the switching conditions.}
\label{carried}
\end{figure}

To understand the iterates of $\varphi$,
we use an associated combinatorial object
called an {\em invariant train track}.
The theory of train tracks is described in \cite{Thurston}; 
we need here  only elementary ideas from this theory.
A {\em train track} is a 3-valent
graph that is embedded on a surface. The edges, called {\em tracks},
are embedded smoothly and the three tangent directions
at the vertices, called {\em switches}, lining up to give a $C^1$-embedding of the
union of any pair of edges meeting at a vertex.
Train tracks have {\em fibered neighborhoods},
closed neighborhoods filled by fibers. Fibers are
intervals transverse to the edges, much
like the tracks of a mono rail, and there is a projection map
of the fibered neighborhood to the train track. 
A curve is {\em carried} by a train track
if it can be isotoped into the fibered neighborhood
so that it is transverse to the fibers.
Such a curve is roughly parallel to the tracks, but may 
run many times over each track.
The curve is determined up to isotopy by a set of {\em weights}.
These are non-negative integers assigned to each track,
giving the number of times the curve
runs over that edge, in either direction.
At each switch there is a {\em switching condition:}
the weight assigned to the one
``incoming'' track is the sum of
the weights of the two ``outgoing'' tracks.
The weights for any two
tracks near a vertex determine the weight for the third. 
An example is shown in Figure~\ref{carried}.

A curve $C$ carried by a train track can be {\em projected}
onto the train track, meaning that the embedding of
the curve can be composed with the projection of each fiber in
the fibered neighborhood to the base point  of that fiber
on the train track.
Each track is given a
weight by the projection, corresponding to the number of pre-images in
$C$ of a point in the interior of the track.
The curve $C$ can be reconstructed from these weights, by taking a
number of copies of each track given by the weights and joining
them together near the switches. There is a unique way to join
that gives an embedded curve. The resulting simple 
closed curve is unique up to isotopy.

\begin{figure}[hbtp]
\centering
\includegraphics[width=.8\textwidth]{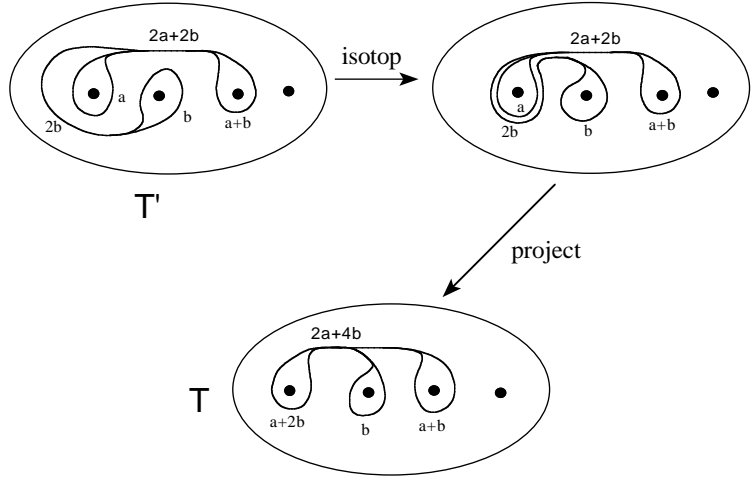}
\caption{A train track $T'$ is projected into another train
track $T'$ that carries it. The weights of the tracks on $T'$ that
project to a given track on $T$ are summed to
give weights on that track.}
\label{project}
\end{figure}

As with curves, a train track $T'$ is
{\em carried} by another train track $T$ if $T'$ can be isotoped
into a fibered neighborhood of $T$ so
that its vertices are carried to vertices
and so that the tracks of $T'$ are transverse to the fibers of the fibered
neighborhood of $T$.
We can then {\em project} $T'$ into $T$
by mapping each fiber to its base point on $T$.
If $T'$ carries weights on its branches, then
these can be summed to give weights on the branches of $T$ to which it
projects, as in Figure~\ref{project}.

A train track is said to be
{\em invariant} under a diffeomorphism $\varphi$ of a surface if its image
$\varphi (T)$ is carried by $T$.

\begin{figure}[hbtp]
\centering
\includegraphics[width=.4\textwidth]{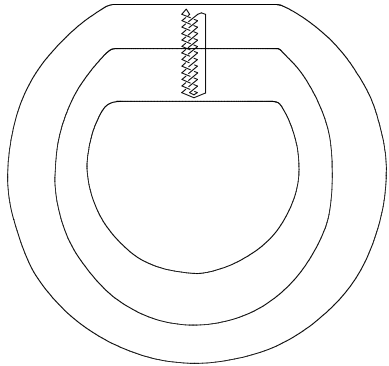}
\caption{Level sets of the function $f_n$ are spheres, but in a ball
containing $K_n$ the level sets are the same as those of the $z$ coordinate.}
\label{St}
\end{figure}

For later application,
we replace the  level planes  $\{z=c\}$ of the height function
$z$ with the level sets of a different function $f_n : \RR^3 \to \RR$, that
agrees with $z$ in a large ball around the origin, a ball that contains
the disks we will be considering.
Thus in subsequent arguments we will be able to
view either $f_n$ or $z$ interchangeably as the
Morse function we are using.
The level sets of $f_n$ are a family of spheres rather than planes $\{z=c\}$.
To construct $f_n$, we first choose a large constant $R_n>0$ such that
a ball of radius $R_n$ centered at the origin contains $F_n$ in its interior.
For each $t$ with $-R_n < t < R_n$, define 
$\Sigma_t$ to be the 2-sphere obtained
by taking the disk $\{(x,y,z): x^2 + y^2 \le R_n, \ z=t\}$
and capping it to form a convex 2-sphere enclosing the point $(0,0,-2R_n)$.
Figure~\ref{St} shows some of these spheres.
The spheres are the level sets of a function 
$$
f_n: \\R^3 \to [-2R_n, \infty).
$$
The restriction of $f_n$ to the disks we will consider
agrees with $z$, 
and $R_n$ will be chosen large enough so that the level sets  of $f_n$
look identical to flat planes in a ball containing the disks.
Note that we can use different functions $f_n$ for different values of $n$,
if necessary, to ensure that our choice of $R_n$ is sufficiently large.
The diffeomorphism of the 4-punctured plane $\varphi$, which
was the identity outside of a disk of radius three around the origin, induces
a diffeomorphism $\varphi : S \to S$ of the
4-punctured sphere $S$, which we call by the same name.

\begin{figure}[hbtp]
\centering
\includegraphics[width=.8\textwidth]{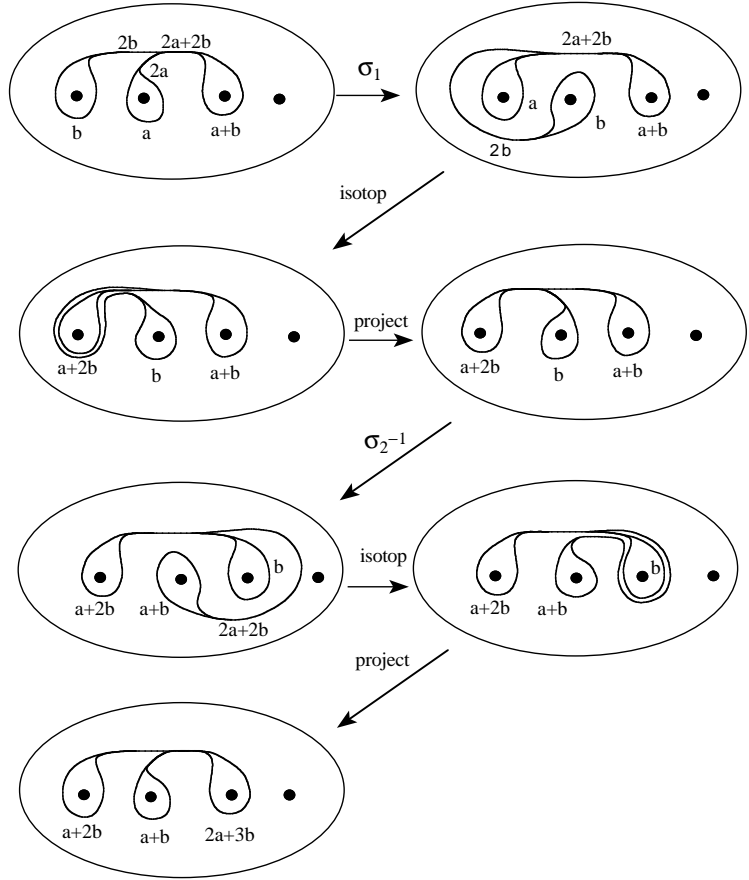}
\caption{An invariant train track $T$ for $\varphi$. Here $\sigma_1$ and
$\sigma_2^{-1}$ denote diffeomorphisms of the four punctured disks
associated to the braid group elements $\sigma_1$ and
$\sigma_2^{-1}$. The composition of these two diffeomorphisms gives
$\varphi$, showing that $T$ is invariant under $\varphi$.}
\label{tt}
\end{figure}

There is an invariant train track $T$ for $\varphi$,
depicted  in  Figure~\ref{tt}, and also shown with
a fibered neighborhood in Figure~\ref{carried}.
An assignment of weights to all the tracks of $T$
is completely determined by assigning
two weights $a$ and $b$
on the two indicated tracks, as in Figure~\ref{carried}.
The non-negative integers $a$ and $b$ are arbitrary,
but all other weights
are determined by the switching conditions.
Each choice of $a, b$ gives rise to a unique simple closed
curve carried by $T$, and
we refer to $a, b$
as the weights with which this curve is carried by $T$.

To understand the iterates of $\varphi$ we study
the image of $T$ under $\varphi$.
The image $\varphi (T)$ can be isotoped so that
vertices of $\varphi (T)$ go to vertices  of $T$ and tracks of
$\varphi (T)$ are transverse to the fibers of the fibered neighborhood
of $T$,  as indicated in Figure~\ref{tt}.

\begin{lemma}
\label{invariant}
The train track $T$ is invariant under the homeomorphism $\varphi$.
A curve carried by $T$ with weights $a, b$ is mapped by
$\varphi$ to a curve carried by $T$
with weights $a+b$ and $a+2b$.
\end{lemma}
{\bf Proof:}
The image of $T$ under $\varphi$ can be isotoped into the fibered
neighborhood of $T$ as shown in Figure~\ref{tt}. The tracks
with initial weights $a$ and $b$ have projected onto
them tracks with total weight $a+b$ and $a + 2b$ respectively.
A curve carried by $T$ with weights $a,b$ is similarly carried
to a curve carried with weights $a+b, a + 2b$.
$~~~\qed$

So $T$ is an invariant train track for $\varphi$, and 
a curve $C$ carried by $T$ with weights
$a$ and $b$, has image
$\varphi(C)$ which is also carried by $T$, but with weights
$a+b$ and $a+2b$.
When $\varphi$ is iterated,
the weights on these two tracks grow according to a Fibonacci sequence:
$$
\{(a,b),( a+b, a+2b),(2a + 3b, 3a + 5b), (5a + 8b, 8a + 11b), \dots \} .
$$

\begin{lemma}
\label{iterated.wts}
A curve carried by $T$ with weights $a_0 \ge 0$
and $ b_0 \ge a_0$ is mapped by the diffeomorphism
$\varphi^n$ to a curve carried by $T$
with weights $a_n$ and $b_n$,
satisfying $a_n \ge 2^n a_0$ and $b_n \ge 2^n b_0$.
\end{lemma}
{\bf Proof:} Under the action of $\varphi$
the weight $a$ corresponding to a curve $C$ is
transformed to the weight $a+b \ge 2a$
corresponding to $\varphi (C)$
and the weight $b$ to $a+2b \ge  2b$.
The result follows by iterating $n$ times. $~~~\qed$ 

\begin{figure}[hbtp]
\centering
\includegraphics[width=.6\textwidth]{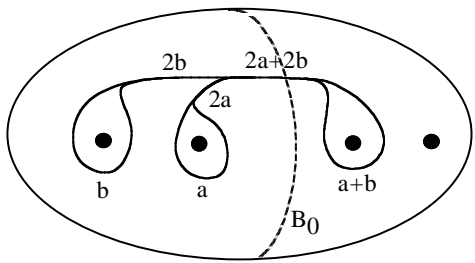}
\caption{The curve $B$ intersects a curve carried by the train track
$T$ in at least $2a+2b$ points.}
\label{B}
\end{figure}

Let $B$ denote the simple closed curve 
on a 4-punctured sphere $S$ that separates
the points $p_1, p_2$ from $p_3, p_4$,
as shown in Figure~\ref{B}.
We analyze the number of intersections between $B$ and a
curve $C$ in the 4-punctured sphere
$S$ that is carried by $T$ with weights $a,b$.
We show there is no isotopy of $C$ in the 4-punctured sphere
which can reduce the number of intersections with $B$ below
$2a + 2b$.

\begin{lemma}
\label{minimal.int}
A curve $C$ in $S$ that is  carried by the train track
$T$ with weights $a$ and $b$ intersects $B$ in at least $2a + 2b$ points.
\end{lemma}
{\bf Proof:} 
In a surface containing two intersecting simple closed curves,
a {\em 2-gon} is a disk on the surface whose boundary consists of an arc from
each of the curves and whose interior is disjoint from each of them.
It is shown in \cite[Lemma 3.1, pp. 108]{HS:85} that if two
simple closed curves on a surface
have more intersections than the minimal possible number
in their isotopy class, then each contains an arc such that
the two arcs together bound a 2-gon on the surface.

Let $C$ be a curve lying in the fibered neighborhood of $T$,
transverse to the fibers and carried with weights $a$ and $b$.
It follows from the above that if
$C$ can be isotoped in $S$ to have fewer
than $2a + 2b$ points of intersection with $B$, then there
exists an arc $\beta$ contained in $B$ and an arc $\gamma$
contained in $C$ that together bound a 2-gon in $S$,
whose interior is disjoint from $B \cup C$. 
We will show that there is no such 2-gon between $C$ and $B$, and 
hence that $C$ can not be isotoped to reduce
the number of its intersections with $B$.

The arc $\gamma$ lies on $C$ and so
lies in the fibered neighborhood of $T$, and
is transverse to the fibers.  Moreover $\gamma$ intersects $B$ only at its
two endpoints, and therefore lies either to the right or to the left
of $B$ on $S$, where ``left'' refers to the side
containing $p_1, p_2$ and ``right'' to the side containing $p_3, p_4$.

An arc carried by $T$ with interior to the right
of $B$ runs once around the third puncture and, together with $\beta$, must
separate the third and fourth punctures.
Similarly an arc carried by $T$ with interior on the left of $B$
runs once around either the first or second
punctures before returning to $B$, and together with $\beta$
separates the first and second punctures.
In either case such an arc is not
homotopic to an arc in $B$, (rel boundary),
and therefore cannot cobound a disk
with an arc $\beta$ contained in $B$.
So $\beta \cup \gamma$ cannot cobound a 2-gon, and it follows that 
the number of intersections of $B$ and $C$ cannot be reduced.
$~~~\qed$

\begin{corollary}
\label{intersection}
Let $p_1, p_2, p_3, p_4$ denote
four distinct marked points on a 2-sphere and let $B$ denote a simple
closed curve separating $p_1, p_2$ from $p_3, p_4$. Let
$\delta$ be the simple closed curve that is the boundary of
a neighborhood of an arc joining $p_1$ to $p_2$ in the complement of $B$.
Then $ \varphi^n (\delta)$ intersects
$B$ in at least $2^{n}$ points.
\end{corollary}
{\bf Proof:}
While $\delta$ is not carried by $T$, its image
$\varphi(\delta)$ is carried by $T$ with weights $a=0, b=1$,
and $\varphi^2(\delta) $ is carried with weights $a=1, b=2$.
Lemma~\ref{iterated.wts} can be
applied to $\varphi^2(\delta) $ and its iterates,
so $\varphi^n (\delta) = \varphi^{n-2} \varphi^2(\delta) $
is carried with weights at least $a= 2^{n-2}, b = 2^{n-1}$.
By Lemma~\ref{minimal.int}, the curve $B$
intersects a curve carried by the train track with weights $a,b$
in at least $2a + 2b$ points.
Since $2a + 2b \ge  2b \ge   2^{n} $, the result follows. $~~~\qed$

\begin{figure}[hbtp]
\centering
\includegraphics[width=.8\textwidth]{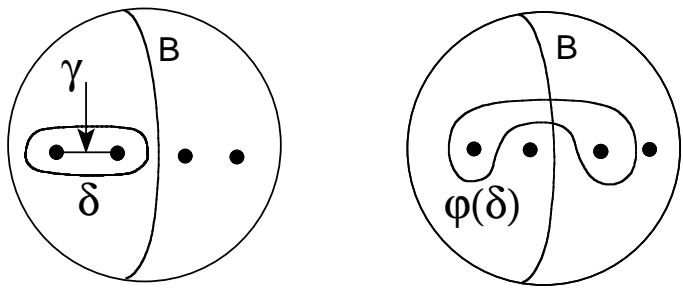}
\caption{The curves $\gamma$, $\delta$ and $\varphi (\delta)$
on a sphere with four marked points removed.}
\label{delta}
\end{figure}

\begin{corollary}
\label{arc.intersection} An arc $\gamma$ joining $p_1$ to $p_2$ in the
complement of $B$ has
image under $\varphi^n$ that intersects the closed curve $B$
on the 4-punctured 2-sphere $S$ in at least $2^{n-1}$ points.
\end{corollary}
{\bf Proof:}
The simple closed curve $\delta$
is isotopic to the boundary of a regular neighborhood
of $\gamma$ and  $\varphi^n(\delta)$ is isotopic to
the boundary of a regular neighborhood of  $\varphi^n(\gamma)$.
For any arc in $S$ that intersects $B$ transversely,
the boundary of a sufficiently thin neighborhood of the arc intersects $B$
in twice the number of points that the arc intersects $B$.
If  $\varphi^n(\gamma)$ intersected $B$ in fewer than $2^{n-1}$ points  
then we can form a thin neighborhood of $\varphi^n(\gamma)$ whose
boundary intersects $B$ in less than $2^n$ points,
contradicting Corollary~\ref{intersection}.
$~~~\qed$

\section{Combinatorial complexity of spanning disks for $K_n$}
\label{comb}

In this section we show that
any PL spanning disk for $K_n$
contains exponentially many triangles, proving
the main result.

{\bf Proof of Theorem~\ref{exp.disks}:}
Let $n$ be any fixed positive integer.
The assertion (1) that $K_n$ is unknotted follows
from its construction as
the composition of a braid and its inverse.

The curve $K_n$ can be constructed with
straight segments as follows:
Four segments above $\{ z=n \}$  and four below $\{ z = -n \}$
cap off the braid.
Between  $\{ z=n \}$  and  $\{ z = -n \}$, a
single line segment forms the entire fourth strand, and the first
three strands are formed from $2n$ copies of the first three strands of $\alpha$.
Each copy  of $\alpha$ requires five
segments for the first three strands.
The total number of segments needed is no more than
$10n + 9 $, which is assertion (2).

We prove assertion (3) in three steps.
Recall that for each fixed $n$,
$K_n$ bounds a smooth disk $F_n$ in $\RR^3$
that we call the standard disk,
and that $B_0 = \Sigma_0 \cap \{ x=0 \}$
is the closed curve obtained by intersecting
$\Sigma_0$ with the plane $ \{x=0\}$.  
We first show that $B_0$
intersects $F_n$ in at least $2^{n-1}$ points.
We then consider an arbitrary smooth spanning disk $E_n$, and show
that the number of intersections of $E_n$
with $B_0$ is at least as large as that of $F_n$
with $B_0$.
In the third step,
we approximate an arbitrary PL disk $D_n$ by
a smooth disk to obtain the same conclusion in the PL setting.

The standard disk $F_n$
is swept out by arcs joining points of $K_n$
in the level sets $\{z = \mbox{ constant} \}$,
as they descend from $\{z>n+2\}$ to $\{z<-n-2\}$.
One arc appears below $\{z=n+2\}$, the second below $\{z=n+1\}$.
The arcs join together to form a single arc at $\{z = -n - 1\}$,
and this in turn disappears below  $\{z = -n - 2\}$.
The height function given by the restriction of the 
$z$-coordinate to $F_n$ defines a Morse function on $F_n$,
and this Morse function has no critical points in the interior of $F_n$.

The arc $\gamma = \gamma_n$ in $F_n \cap \{z=n\}$
connects $p_1$ and $p_2$, as shown in 
Figure~\ref{delta}.
Denote by $\gamma_t$ the arc in $F_n \cap \{z=t\}$ that is
in the same component of $F_n \cap \{z \ge t\}$ as $\gamma$.
For each integer $k$ with $0<k \le n$, as
we slide  $\gamma_k$  down  one unit along $F_n$,
$\gamma_{k} $ is deformed along $F_n$ to
$ \gamma_{k-1} = \varphi ( \gamma_k )$. So as $t$ decreases to 0,
the arc $ \gamma_n $ is slid along $F_n$
to an arc $\gamma_0$ that
is the image of $n$ iterations of $\varphi$.

Let $B_0$ denote the closed curve
along which the level set $\Sigma_0 = {f_n}^{-1}(0)$ intersects the $yz$-plane.
Then $B_0$ separates the four points of intersection of $\Sigma_0$
and $K_n$ into pairs, $p_1,p_2$ and $p_3,p_4$.
The standard disk $F_n$ intersects $B_0$ in at least $2^{n-1}$ points by
Corollary~\ref{arc.intersection}.

Our goal is to show that an arbitrary PL disk
bounded by $K_n$ intersects $B_0$ in at least as many points as does $F_n$.
Before considering PL disks, we first consider a
smooth spanning disk $E_n$.
In this setting we will apply some basic
results from the Morse Theory of smooth functions on
surfaces;
see \cite{Milnor} for an exposition of smooth Morse Theory.
We will then shift back to the PL setting.
The height function $z$, or the function $f_n$ that agrees with
it on the region we are studying, will serve as the Morse functions.
Let $E_n$ denote an arbitrary disk spanning $K_n$ such that
\begin{enumerate}
\item $E_n$ has smoothly embedded interior.
\item The height function $z$ restricted to $E_n$ is a Morse function.
\end{enumerate}

We now show, using Morse theory, that the surface $E_n$ 
intersects the closed curve $B_0$ in at least as many
points as does the ``standard disk'' $F_n$.
Choose a value of $R_n$ large enough 
so that $F_n$ and $E_n$
both lie in the interior of the ball of radius $R_n$, and
as before form the Morse function $f_n$
whose level sets are spheres $\Sigma_t$ for $t > -2R_n$.
The intersection of $E_n$ with the spheres $\Sigma_t$ at
non-critical levels is contained in $\Sigma_t \cap \{z=t\}$.
As $t$ decreases from $\infty$ to $-2R_n$, 
the sphere $\Sigma_t$ begins to intersect $K_n = \partial E_n$,
when $t = n+2$.
As $t$ decreases there are first one,
then two arcs in $\Sigma_t \cap E_n$, along with a (possibly
empty) collection of simple closed curves. 
For $n+1 < t < n+2 $, $\Sigma_t \cap E_n$ consists of a
single arc $\beta_t$, along with a
possibly empty collection of simple closed curves.
For $ -n < t <  n$, $\Sigma_t \cap E_n$ contains two arcs connecting the four points
of $\Sigma_t \cap K_n$.
As $t$ decreases from $n+1$, $\beta_t$ is 
continuously deformed as long as $E_n$ is transverse to $\Sigma_t$.
As long as the transversality continues to hold, let
$\beta_t$ denote the arc that is in the same component
of $E_n \cap \{z \ge t\}$ as $\beta$.

As long as passing through the critical level \{$t=c$\}
does not change which pairs of points on $K_n$ are connected by
the pairs of arcs,
it is possible to extend the definition of
$\beta_t$ to one of the two arcs below the critical point.
We define $\beta_{c-\epsilon}$ to be the arc connecting the
same pair of points as  $\beta_{c+\epsilon}$.
In these cases even the isotopy class of the arc is preserved, though
the curve $\beta_t$ does not change continuously
when $t$ passes through the critical level $c$.

\begin{figure}[hbtp]
\centering 
\includegraphics[width=.75\textwidth]{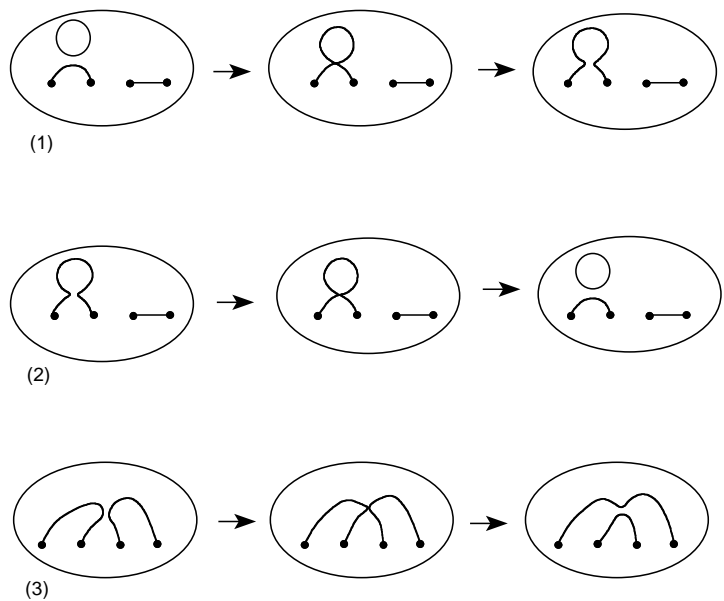}
\caption{Three types of changes  can occur
in the level sets $E_n \cap \{z = t\}$ as they pass through a critical point.}
\label{arcchanges}
\end{figure}

For the values $-n < t < n$, the level sets of $f_n$ 
are transverse to $K_n$, so
any critical points lie in the interior of
the disk $E_n$.
There are three types of changes in $\beta_t$ that can
occur when descending from
$\Sigma_{c+\epsilon}$ to $\Sigma_{c-\epsilon}$, as indicated in Figure~\ref{arcchanges}.

\begin{enumerate}
\item
Moving past a saddle critical point connects
$\beta_{c+\epsilon}$ to a simple closed curve of $\Sigma_{c+\epsilon} \cap E_n$ 
to form $\beta_{c- \epsilon}$.
\item
Moving past a saddle critical point connects $\beta_{c+\epsilon}$ to itself
to form $\beta_{c- \epsilon}$ together with a simple closed curve.
\item
Moving past a saddle critical point connects $\beta_{c+\epsilon}$ to
the second arc of $\Sigma_{c+\epsilon} \cap E_n$.
No arc $\beta_{c- \epsilon}$ is defined.
\end{enumerate}

The first two types of moves are inverses, since reversing
the direction of a type (1) move gives a type (2) move and vice-versa. 
The level surface in which $\beta_{c+\epsilon}$ lies is either a
2-punctured sphere,
if $c \ge n+1$ or $c \le -n-1$, and is a 4-punctured sphere otherwise.
In a 2-punctured sphere there is a unique isotopy class of
arcs connecting the two punctures, so the
isotopy class of $\beta_t$ in $\Sigma_t$ is unchanged when passing through the
critical point.
In a 4-punctured sphere there
are many isotopy classes of arcs connecting two of the punctures,
but as we pass a saddle critical point of type (1) or type (2),
the curve $\beta_t$
remains in the complement of the second arc of intersection.
The complement of an arc in a sphere
is homeomorphic to an (open) disk,
and in a disk there is a unique isotopy class of arcs
connecting any two points.
So the isotopy class of $\beta_t$
in  $\Sigma_t$ remains unchanged by a saddle move in these cases. 

The third type of critical point
does change the isotopy class of the arc, since it
changes the boundary points connected
by the arc. The following lemma asserts that this
can occur at most once.

\begin{figure}[hbtp]
\centering
\includegraphics[width=.5\textwidth]{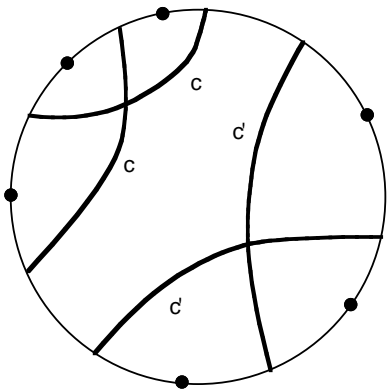}
\caption{The existence of two type (3) saddle critical points
on a disk, at levels ${f_n}^{-1}(c)$ and ${f_n}^{-1}(c')$,
imply that the boundary of the disk
has at least six critical points.}
\label{2type3}
\end{figure}

\begin{lemma}
\label{2crit}
Suppose that $f_n:D \to \RR$
is a  Morse function on a topological disk $D$ that
 restricts to a Morse function on $\partial D$.
Suppose also that  $f_n|_{\partial D}$, the restriction of $f_n$ to $\partial D$,
has at most four critical points on $\partial D$.
Then $f_n$ can have at most one interior critical point of
type (3) that is a saddle 
connecting distinct arcs in the level set of $f_n$. 
\end{lemma}
{\bf Proof:} 
Suppose there is an interior critical point with critical value $c$
that is a saddle connecting distinct arcs in the level set $f_n=c+ \epsilon$.
The four arcs leaving the saddle point hit the boundary of
$D$ at four distinct points. These arcs divide $D$ into four
quadrants, which meet $\partial D$ in four arcs. Each of the
boundary arcs has its two endpoints on the level set $f_n=c$,
and is non-constant on these boundary arcs.
So each contains at least one maximum or minimum of $f_n|_{\partial D}$.
Suppose there were a second saddle critical point
on a level set $\{ f_n = c' \}$.
Since critical points of Morse functions have
distinct values, we have that $c' \ne c$.
The four arcs emerging from the second saddle are
therefore disjoint from the first critical level set, and
contained in one of the previously defined quadrants. 
See Figure~\ref{2type3}.

The intersection of this quadrant with $\partial D$ has four points
on which $f_n$ takes the value $c'$, and therefore $f_n|_{\partial D}$
has at least three critical points in this quadrant.  It follows that
$f_n|_{\partial D}$ has at least six critical points, contradicting
the hypothesis.  So only one saddle of type (3) can occur.
$~~~\qed$

Since $E_n$ is a topological disk whose boundary has four
critical points for the height function,
at most one type (3) critical point can occur.
Assume first that a type (3) critical point does not occur for $t > 0$, 
On the standard disk, $\gamma_0$ is obtained
from $\gamma_n$ by a continuous
deformation involving no critical points, while $\beta_0$
is obtained from $\beta_n$ by a process that may
include passing through 
critical points of types (1) and (2), but none
of type (3).
Therefore the isotopy class of $\beta_0$ in the
4-punctured sphere is the same as that of $\beta_{n+2-\epsilon}$.
Since $\beta_t$ is isotopic to $\gamma_t$ for $t$ close to $n+2$, we conclude
that $\beta_0$ is isotopic to $\gamma_0$.
By Lemma~\ref{minimal.int} $\beta_0$ intersects $B_0$
in at least as many points as $\gamma_0$.

Now consider the case where a type (3) critical
point does occur for some $t > 0$. By Lemma~\ref{2crit},
there are no type (3) critical points for $ t <0$.
In this case we repeat the previous argument,but using
the function $-f_n$ rather than $f_n$.
Note that as
$z$ increases from $\{ z = k \}$
to  $\{ z = k+1 \}$ for $k$ an integer with $-n \le k \le -1$,
the level sets of $F_n$ are again transformed by an application of
$\varphi$.
We replace $\gamma$ with the
arc $\gamma'$ in $E_n \cap \Sigma_{-n}$ which joins $p_2$ to $p_3$,
and $\delta$ with $\delta'$,
the boundary of a regular neighborhood of $\gamma_n'$.

\begin{figure}[hbtp]
\centering
\includegraphics[width=.6\textwidth]{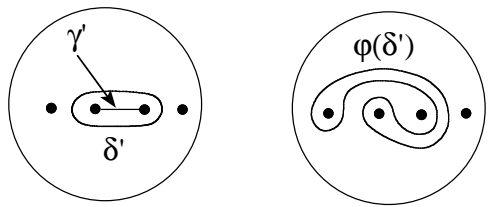}
\caption{The curves $\gamma', \ \delta'$ and $\varphi (\delta')$.}
\label{delta'}
\end{figure}

The curve $\delta'$ is carried with weights $a=1, b=0$ by $T$
and $\varphi (\delta') $ is carried with weights
$a=1, b=1$. By Lemma~\ref{iterated.wts}, $\varphi^n (\delta')$
is carried by $T$ with weights $a \ge 2^{n-1}$ and
$b \ge 2^{n-1}$. 
By Lemma~\ref{minimal.int}
a curve isotopic to  $\varphi^n (\delta')$ intersects $B_0$
in at least $2a + 2b = 2^{n+1}$ points.
Then  $\varphi^n (\gamma')$ intersects $B_0$ in at least half as many,
or $2^n$ points.
As $t$ increases from $-n$ to 0, 
the arc $\gamma_t'$ in the component of $E_n \cap \{ z \le t\}$
containing $\gamma'$, is carried to
an arc in $\Sigma_0$ isotopic to $\varphi^n (\gamma_n')$.
So $E_n$ again must intersect $B_0$ in at least $2^{n-1}$ points.
In every case the curve $B_0$ intersects $E_n$
in at least $2^{n-1}$ points.

Now consider an arbitrary PL disk $D_n$ with boundary $K_n$.
After an arbitrarily small isometry of $\RR^3$,
we can arrange that $D_n$ intersects the $y$-axis transversely in
a finite number of points,
The disk $D_n$ can be approximated by a
disk $E_n$, with smoothly embedded interior, that
coincides with $D_n$ in a neighborhood of each intersection point with the
$y$-axis, and that remains disjoint from other points of 
the $y$-axis.
We choose $R_n$ larger, if necessary, so that in the passage from
the Morse function $z$ to the Morse function
$f_n$, the three surfaces $D_n, E_n$ and $F_n$ that we consider
only intersect the flat parts of the spheres $S_t$ which form
the level sets of $f_n$.
So the intersection of $B_0$ with $D_n, E_n$ and $F_n$ is the
same as that of the $y$-axis with these surfaces,
each intersection being in the interior of a triangular face of $D_n$.

But we have shown that $E_n$ intersects 
$B_0$ in at least $2^{n-1}$ points, and it therefore follows that
the PL disk $D_n$ also intersects $B_0$ in at least $2^{n-1}$ points. 

Since a triangle transversely intersects a line in at most one point,
and $B_0$ agrees with  the $y$-axis
in a ball containing $D_n$, $E_n$ and $F_n$,
this implies that $D_n$ contains at least $2^{n-1}$ triangles,
and Theorem~\ref{exp.disks} is proved.
$~~~~~ \qed$

\noindent
{\bf Remarks:}
\begin{enumerate}

\item If we allow spanning disks that self-intersect,
then the number of triangles required to span $K_n$ grows only linearly with $n$.
If a spanning surface of arbitrary genus is allowed,
it can be shown that the number of triangles required to
span $K_n$ grows at most quadratically in $n$ \cite{HL:01}.

\item The number of Reidemeister moves required to
transform any unknotted curve constructed with $n$ polygonal edges
into a single triangle
has an exponential upper bound derived in \cite{HL}.
For the particular $K_n$ constructed here,
the number of Reidemeister moves required to transform
the projection of
$K_n$ to a projection with no crossings grows only linearly with $n$.

\item The argument establishes somewhat better estimates then claimed above.
If we embed $K_n$ into $\RR^3$ in a more efficient way, then for large $n$,
$K_n$ has at most $6n$ segments and the number of triangles
contained in any disk spanning
$K_n$ grows faster than a constant times $\phi^{2n}$,
where $\phi$ is the golden ratio. 

\end{enumerate}

\section{Acknowledgments}
This paper was completed while the first author was visiting
the Institute for Advanced Study.
The authors are grateful to J. Lagarias and the referee for
helpful suggestions on the exposition.

\end{document}